
\magnification=1200

\font\tengoth=eufm10
\font\sevengoth=eufm7
\newfam\gothfam
\textfont\gothfam=\tengoth
\scriptfont\gothfam=\sevengoth
\def\goth{\fam\gothfam\tengoth}

\font\tenbboard=msbm10
\font\sevenbboard=msbm7
\newfam\bboardfam
\textfont\bboardfam=\tenbboard
\scriptfont\bboardfam=\sevenbboard
\def\bboard{\fam\bboardfam\tenbboard}

\newif\ifpagetitre
\newtoks\auteurcourant \auteurcourant={\hfill}
\newtoks\titrecourant \titrecourant={\hfill}

\pretolerance=500 \tolerance=1000 \brokenpenalty=5000
\newdimen\hmargehaute \hmargehaute=0cm
\newdimen\lpage \lpage=14.3cm
\newdimen\hpage \hpage=20cm
\newdimen\lmargeext \lmargeext=1cm
\hsize=11.25cm
\vsize=18cm
\parskip=0cm
\parindent=12pt

\def\margehaute{\vbox to \hmargehaute{\vss}}
\def\margebasse{\vss}

\output{\shipout\vbox to \hpage{\margehaute\nointerlineskip
  \corpsdepage\margebasse}
  \advancepageno \global\pagetitrefalse
  \ifnum\outputpenalty>-2000 \else\dosupereject\fi}

\def\corpsdepage{\hbox to \lpage{\hss\pagetexte\hskip\lmargeext}}
\def\pagetexte{\vbox{\makeheadline\pagebody\makefootline}}
\headline={\ifpagetitre\titleheadline \else
  \ifodd\pageno\rightheadline \else \leftheadline\fi\fi}
\def\leftheadline{\hfil\the\auteurcourant\hfil}
\def\rightheadline{\hfil\the\titrecourant\hfil}
\def\titleheadline{\hfill}
\pagetitretrue

\font\petcap=cmcsc10
\def\pc#1#2|{{\tenrm#1\sevenrm#2}}
\def\pd#1 {{\pc#1|}\ }

\def\pointir{\discretionary{.}{}{.\kern.35em---\kern.7em}\nobreak\hskip 0em
 plus .3em minus .4em }

\def\titre#1|{\message{#1}
             \par\vskip 30pt plus 24pt minus 3pt\penalty -1000
             \vskip 0pt plus -24pt minus 3pt\penalty -1000
             \centerline{\bf #1}
             \vskip 5pt
             \penalty 10000 }
\def\section#1|
               {\par\vskip .3cm
               {\bf #1}\pointir}
\def\ssection#1|
             {\par\vskip .2cm
             {\it #1}\pointir}

\long\def\th#1|#2\finth{\par\vskip 5pt
              {\petcap #1\pointir}{\it #2}\par\vskip 3pt}
\long\def\tha#1|#2\fintha{\par\vskip 5pt
               {\petcap #1.}\par\nobreak{\it #2}\par\vskip 3pt}

\def\rem#1|
{\par\vskip 5pt
                {{\it #1}\pointir}}
\def\rema#1|
{\par\vskip 5pt
                {{\it #1.}\par\nobreak }}

\def\qed{\quad\raise -2pt \hbox{\vrule\vbox to 10pt
{\hrule width 4pt
\vfill\hrule}\vrule}}

\def\cqfd{\ifmmode
\unkern\quad\hfill\qed
\else
\unskip\quad\hfill\qed\bigskip
\fi}

\newcount\n
\def\exo{\advance\n by 1 \par \vskip .3cm {\bf \the \n}. }



\hfill

\vskip 20 pt

\centerline{ANALYSIS OF THE MINIMAL REPRESENTATION  OF ${\rm Sp}(r,{\bboard R})$  }

\vskip 20 pt
 
\centerline{Dehbia Achab}

\vskip 4 pt

\centerline{\tt dehbia.achab@imj-prg.fr}

\vskip 10 pt

\vskip 15pt

\noindent
{\bf Abstract}   The minimal representations of ${\rm Sp}(r,{\bboard R})$ can be realized on a Hilbert space of holomorphic functions. This is the analogue of the Brylinski-Kostant model. It can also be realized on a Hilbert space of $L^2$ functions on ${\bboard R}^n$. This is the Schr\"odinger model. We will describe the two realizations and  a transformation which maps one model to the other. It involves the classical Bargmann transform.


\vskip 4 pt
\noindent
\bf Introduction. \rm   In the paper [A11]   a general construction  for a simple  complex Lie  algebra ${\goth g}$ and of a real form ${\goth g}_{\bboard R}$  has been proposed, starting from a pair $(V,Q)$  where $V$ is a simple Jordan algebra of rank $r$  and $Q$ a polynomial on $V$, homogeneous of degree $2r$.  The Lie algebra ${\goth g}$ is of Hermitian type. 
In the paper [A12b],  the manifolds $\Xi$ and $\Xi^{\sigma}$ are   the orbits of  the linear form $\tau_{c_{\lambda}}$  and its conjugate $\tau_{c_{\lambda}}^{\sigma}=\kappa(\sigma)\tau_{c_{\lambda}}$, for an idempotent $c_{\lambda}$,   under the structure group ${\rm Str}(V)$  acting  by the restriction of an irreducible representation ${\kappa}$ of the conformal group ${\rm Conf}(V,Q)$.    The spaces ${\cal F}(\Xi)$  and ${\cal F}(\Xi^{\sigma})$ are  Hilbert spaces  of holomorphic functions on the complex manifolds $\Xi$ and $\Xi^{\sigma}$. The minimal representations of ${\goth g}_{\bboard R}$ are realized in ${\cal F}(\Xi)$ and ${\cal F}(\Xi^{\sigma})$. 
In this paper we consider the special case where $V={\rm Sym}(r,{\bboard C})$,  $Q$ is the square  of the determinant and the construction leads to the Lie algebra ${\goth g}={\goth sp}(r,{\bboard C})$ and to a real form ${\goth g}_{\bboard R}=Ad(g_0)({\goth sp}(r,{\bboard R}))$ for some $g_0\in {\rm Sp}(r,{\bboard C})$. The two  minimal representations  $\rho$  and $\rho^{\sigma}$ of  ${\goth g}_{\bboard R}$ are    respectively realized   in ${\cal F}(\Xi)$ and ${\cal F}(\Xi^{\sigma})$, which turn out to be the classical Fock space ${\cal F}({\bboard C}^r)$.  They integrate to representations $T$ and $T^{\sigma}$ of the complex group ${\rm Mp}(r,{\bboard C})$ and their restriction to $Mp(r,{\bboard R})$ realize unitary representations in $T(g_0^{-1})({\cal F}({\bboard C}^r))$ and $T^{\sigma}(g_0^{-1})({\cal F}({\bboard C}^r))$. 
On another hand, we consider  the two minimal representations $R$ and $R^{\sigma}$, i.e. the Segal-Shale-Weil representations,  of ${\rm Mp}(r,{\bboard R})$  on  $L^2({\bboard R}^r)$. 
 We   also  give   unitary  integral operators ${\cal B}$   from the space $L^2({\bboard R}^r)$ onto the  space $T(g_0^{-1})({\cal F}(\Xi))$   which intertwins   $\rho=dT$ and $dR$ of ${\goth sp}(r,{\bboard R})$,  and ${\cal B}^{\sigma}$ from the space $L^2({\bboard R}^r)$ onto the  space $T^{\sigma}(g_0^{-1}){\cal F}(\Xi^{\sigma})$ which intertwins $\rho^{\sigma}=dT^{\sigma}$ and $dR^{\sigma}$ of  ${\goth sp}(r,{\bboard R})$. 
\bigskip
\noindent
\bf 1. The analogue of the Brylinski-Kostant model. General case \rm 
\vskip 6 pt
\noindent
Let $V$ be a simple complex Jordan algebra with rank $r$ and  dimension $n$ and $Q$ the 
homogeneous polynomial of degree $2r$ on $V$ given by  $Q(v)=\Delta(v)^2$ where $\Delta$ is the Jordan algebra determinant.  The structure group of $V$ is defined by:
\vskip 4 pt

\centerline{$ {\rm Str}(V)=\{g\in {\rm GL}(V) \mid \exists \chi(g)\in {\bboard C},
\Delta(gz)=\chi(g)\Delta(z)\}$}

\vskip 4 pt
\noindent
The conformal group ${\rm Conf}(V)$ is the group of rational transformations $g$ of $V$ generated by: the translations $z\mapsto z+a$ ($a\in V$), the dilations $z\mapsto \ell z$
($\ell \in {\rm Str}(V)$), and the  conformal inversion $\sigma : z\mapsto -z^{-1}$ (see [M78]). 
\vskip 1 pt
Let ${\goth p}$ be the space of polynomials on $V$ 
generated by the polynomials  $Q(z-a)$ of $Q$, with $a\in V$.  Let $\kappa$ be the cocycle representation of  $K={\rm Conf}(V)$,  defined in [A11] and [AF12] as follows:
\vskip 4 pt

\centerline{
$(\kappa(g)p)(z)=\mu(g^{-1},z)p(g^{-1} z),$}

\vskip 2 pt

\centerline{$\mu(g,z)=\chi((Dg(z)^{-1})\quad (g\in K,\ z\in V).$}

\vskip 4 pt
\noindent
The function $\kappa(g)p$ belongs actually to $\goth p$  (see [FG96], Proposition 6.2).
The cocycle $\mu(g,z)$ is a polynomial in $z$ of degree 
$\leq \, {\rm deg}\, Q$ and
$$\eqalign{
(\kappa(\tau_a)p)(z)&=p(z-a) \quad (a\in V),\cr
(\kappa(\ell )p)(z)&=\chi(\ell )p(\ell^{-1} z) \quad (\ell\in L),\cr
(\kappa(\sigma)p)(z)&=Q(z)p(-z^{-1}).\cr}$$
Let $L={\rm Str}(V)$. 
It is established in [A11] that the element  $H_0\in{\goth z}({\goth l})$,  given by ${\rm exp}(tH_0)=l_{e^{-t}}: z\in V\mapsto e^{-t}z$ , defines  a grading  of ${\goth p}$:
\vskip 2 pt

\centerline{${\goth p}={\goth p}_{-r}\oplus{\goth p}_{-r+1}\oplus\ldots\oplus{\goth p}_0\oplus\ldots\oplus{\goth p}_{r-1}\oplus{\goth p}_{r},$}

\vskip 2 pt
 
\centerline{${\goth p}_j=\{p\in {\goth p} \mid {\rm d}\kappa(H_0)p=jp\}$}

\vskip 4 pt
\noindent
 is the set of polynomials in $\goth p$,  homogeneous  of degree 
$j+r$.  Furthermore 
$\kappa(\sigma) : {\goth p}_j \rightarrow
{\goth p}_{-j}$, and
${\goth p}_{-r}={\bboard C}, \quad {\goth p}_r={\bboard C}\, Q, 
\quad {\goth p}_{r-1} \simeq V,\quad {\goth p}_{-r+1}\simeq V.$ 
Observe that  ${\goth p}_{r-1}=\{\kappa(\sigma)p \mid p\in  {\goth p}_{-r+1}\}$ and that 
$ {\goth p}_{-r+1}$ is the space of linear forms on $V.$
\vskip 4 pt
 Assume $r\ne 1$ and  denote by ${\cal W}={\cal V}\oplus {\cal V}^{\sigma}$ with  
${\cal V}= {\goth p}_{-r+1},  {\cal V}^{\sigma}= {\goth p}_{r-1}.$
\vskip 2 pt
\noindent
Let's consider the linear form   $\tau: z\mapsto {\rm tr}(z)={\rm trace}(z)$ and its image $\tau^{\sigma}$ given by  $\tau^{\sigma}(z)=Q(z)\tau(-z^{-1})$. 
Denote by 
${E}=\tau$ and ${F}=\tau^{\sigma}$ and let $X_{{E}}\in {\goth k}_1$ and  $X_{{F}}\in{\goth k}_{-1}$  such that  
${E}=$d$\kappa(X_{E})1$ and ${F}=$d$\kappa(X_{F})Q$.  Put $[1,Q]=H_0$, $[E,Q]=-X_E$, $[F,1]=X_F$. Let $\lambda_0$ be  a bilinear form  on ${\cal V}\times{\cal V}^{\sigma}$.
Then   ${\goth g}={\goth l}\oplus {\cal  W}$   carries  a unique simple Lie algebra structure such that 
$$\eqalignno{
({\rm i}) \quad &[X,X']=[X,X']_{\goth l} \quad (X,X' \in {\goth l}), \cr
({\rm ii}) \quad &[X,p]={\rm d}\kappa(X)p \quad (X\in {\goth l},  p\in {\cal W}), \cr
({\rm iii})\quad &[{E},{F}]=\lambda_0(E,F)H_0+[X_E,X_F].}$$
\vskip 4 pt
We recall also  the  real form ${\goth g}_{\bboard R}$ of $\goth g$ which will be considered in the sequel.   We fix a  Euclidean real form $V_{\bboard R}$ of the complex Jordan algebra $V$, denote by $z\mapsto
\bar z$ the conjugation of $V$ with respect to $V_{\bboard R}$, and then 
consider  the involution $g\mapsto \bar g$ of $K$ given  by: 
$\bar g z=\overline{g\bar z}$.  
  The involution $\alpha$ defined by $\alpha (g)=\sigma\bar g\sigma^{-1}$ is a Cartan involution of $K$  and  $K_{\bboard R}=\{g\in K\mid \alpha(g)=g\}$
is a compact real form of $K$ and it follows that   $L_{{\bboard R}}=L\cap K_{{\bboard R}}$ is a compact real form of $L$. 
Observe that, since for $g\in {\rm Str}(V)$, $\sigma\circ g\circ \sigma=g'$, the adjoint of $g$ with respect to the symmetric  form $(w\mid w')=\tau(w\bar w')$, then 
$L_{{\bboard R}}=\{l\in L \mid ll'={\rm id}_V\}.$
Let ${\goth u}$ be the  compact real form of $\goth g$ such that 
${\goth l}\cap {\goth u}={\goth l}_{\bboard R}$, the Lie algebra of $L_{\bboard R}$.  Denote by ${\cal W}_{{\bboard R}}={\goth u}\cap(i{\goth u})$. Then, the real Lie algebra  defined by
$\goth g_{\bboard R}={\goth l}_{\bboard R}+{\cal W}_{\bboard R}$
 is a real form of $\goth g$ and this  decomposition  is its Cartan decomposition.
  Since the complexification of the Cartan decomposition of ${\goth g}_{\bboard R}$ is $\goth g={\goth l}+{\cal W}$ and since ${\cal W}={\cal V}+{\cal V}^{\sigma}={\rm d}\kappa({\cal U}({\goth l})){E}+{\rm d}\kappa({\cal U}({\goth l})){F}$ is a sum of two simple ${\goth l}$-modules, it follows that the simple real Lie algebra ${\goth g}_{\bboard  R}$ is of  Hermitian type. One can  show that 
${\cal W}_{\bboard R}=\{p\in {\cal W}\mid \beta (p)=p\}$ 
where  we defined for a polynomial  $p\in {\cal W}$, $\bar p=\overline{p(\bar z)},$ and 
 considered the antilinear involution $\beta $ of ${\cal W}$ given by $\beta (p)=\kappa (\sigma )\bar p.$ 
 \vskip 4 pt
 Consider  the isomorphisms   $v\in V\mapsto \tau_v\in {\cal V}$ and $v\in V\mapsto \tau_v^{\sigma}\in {\cal V}^{\sigma}$ where  $\tau_v$ is the linear form on $V$ given by $\tau_v(v')=\tau(vv')$ and $\tau_v^{\sigma}=\kappa(\sigma)\tau_v$.
For a suitable choice of a bilinear form  $\lambda_0: {\cal V}\times{\cal V}^{\sigma} \rightarrow {\bboard C}$,  ${\goth g}$ and ${\goth g}_{\bboard R}$ are   isomorphic to  matrix Lie algebras ([A12b]).
  Let ${\goth h}$ be a Cartan subalgebra of ${\goth l}$, which contains $H_0$. Since the Lie algebra ${\goth g}={\goth l}+{\cal W}$ is of Hermitian type, then ${\goth h}$ is also  a Cartan subalgebra of ${\goth g}$.  Let  $\Delta^+({\goth g},{{\goth h}}) $ be  a suitably chosen positive system of roots and let   $\lambda$ be the highest root. Denote by ${\goth g}_{\lambda}$  and  ${\goth g}_{-\lambda}$ the root spaces  corresponding to the highest  and  lowest root. 
One knows that the minimal adjoint nilpotent orbit of ${\goth g}$ is given by ${{\cal O}}_{{\rm min}}=G\cdot ({\goth g}_{\lambda}\backslash\{0\})$, where $G={\rm Int}({\goth g})$ is the group of inner automorphisms of ${\goth g}$.   Recall from [A12b], Lemma 2.1, that ${\goth g}_{\lambda}\subset {{\cal V}}\backslash\{0\}$ and  ${\goth g}_{-\lambda}=\{{\kappa(\sigma)T} \mid  T\in{\goth g}_{\lambda}\} \subset {\cal V}^{\sigma}\backslash\{0\}.$
Denote by $\tau_{\lambda}$ the linear form given by $\tau_{\lambda}(v)=\tau_{c_{\lambda}}(v)=\tau(c_{\lambda}v),$ and  by $\tau_{\lambda}^{\sigma}=\kappa(\sigma)\tau_{\lambda}$. 
  They are 
nilpotent elements of ${\goth g}$, $\tau_{\lambda} \in {\goth g}_{\lambda}$ and  $\tau_{\lambda}^{\sigma} \in {\goth g}_{-\lambda}.$
The orbits $\Xi$ of $\tau_{\lambda}$ and $\Xi^{\sigma}$ of $\tau_{\lambda}^{\sigma}$ under the group $L$   acting on ${\goth p}$ by the representation $\kappa$,   are conical varieties related by $\Xi^{\sigma}=\kappa(\sigma)\Xi$.  They are  the minimal nilpotent $L$-orbits  in ${\cal W}$.
\vskip 1 pt
\noindent
Polynomials $ \xi\in \Xi$   and   $ \xi^{\sigma}=\kappa(\sigma)\xi\in \Xi^{\sigma}$ can be written 
\vskip 4 pt

\centerline{$\xi(v)=\kappa(l)\tau_{\lambda}(v)=\chi(l)\tau(((l')^{-1}c_{\lambda})v),$}
 
 \vskip 1 pt
 \noindent
 
\centerline{$\xi^{\sigma}(v)=\kappa(\sigma)\xi(v)=\chi(l)Q(v)\tau(-((l')^{-1}c_{\lambda})v^{-1}),$}

\vskip 1 pt
\noindent
where $l'\in{\rm Str}(V)$ is   the adjoint of $l$ for  the  inner product  $(x\mid y)={\rm tr}(xy)$. 
Then  $\Xi$ and $\Xi^{\sigma}$ are   realized  as a $L$-orbit of $c_{\lambda}$ in $V$ and have  explicit coordinate systems.
\vskip 4 pt
\hfill
\eject
\noindent
We consider in this paper the  case of 
 $V_{\bboard R}={\rm Sym}(r,{\bboard R})$,  $V={\rm Sym}(r,{\bboard C})$,      ${\rm Str}(V_{\bboard R})={\rm GL}(r,{\bboard R})$,   ${\rm Str}(V)={\rm GL}(r,{\bboard C})$ acting on $V_{\bboard R}$ or on $V$ by $gx=g\cdot x\cdot g^t$. Then $c_{\lambda}={\rm diag}(1,0,\ldots,0)$,   the orbits  $\Xi$ and $\Xi^{\sigma}$ are both   identified to 
$$\Gamma(\Xi)=\{g\cdot c_{\lambda}\cdot g^t \mid g\in {\rm GL}(r,{\bboard C})\}=\{\xi_z=zz^t \mid z\in {\bboard C}^r\}.$$
 
\vskip 4 pt
\noindent
The group $L$  acts on the spaces ${{\cal O}}(\Xi)$  and  ${{\cal O}}(\Xi^{\sigma})$ of holomorphic functions on $\Xi$ and $\Xi^{\sigma}$ respectively
by:

\vskip 1pt

\centerline{$\bigl(\pi_{\alpha}(l)f\bigr)(\xi)=\chi(l)^{\alpha}f\bigl(\kappa(l)^{-1}\xi\bigr) $ and $\bigl(\pi_{\alpha}^{\sigma}(l)f\bigr)(\xi^{\sigma})= \chi(l)^{\alpha}f\bigl(\kappa(l)^{-1}\xi^{\sigma}\bigr).$}

 \vskip 2 pt
 \noindent
For  $\xi=X_z\in \Gamma(\Xi)$  and for every function $f\in {{\cal O}}(\Xi)$,  we write $f(\xi)=\phi(z)$. In these coordinates, the representations $\pi_{\alpha}$ and $\pi_{\alpha}^{\sigma}$   are   given by

\vskip 4 pt

\centerline{$\pi_{\alpha}(l) \phi(z)= \chi(l)^{\alpha}\phi(l_1^{-1}\cdot z)$ for $l=(l_1,(l_1^t)^{-1})$}

\vskip 2 pt
\noindent
 and 
 \vskip 2 pt
 
 \centerline{$\pi_{\alpha}^{\sigma}(l) \phi(z) =\chi(l)^{-\alpha}\phi(l_1^{t}\cdot z)$ for  $l=(l_1,(l_1^t)^{-1})$.}

\vskip 2 pt
\noindent

\vskip 1 pt
\noindent

  \vskip 2 pt
For    $m\in {\bboard Z}$, let 
${\cal O}_m(\Xi)$ and     ${\cal O}_m(\Xi^{\sigma})$  be the the spaces of holomorphic functions $f$ on   $\Xi$  and $f^{\sigma}$ on  $\Xi^{\sigma}$ respectively  such that  for every $w\in {\bboard C}^*$,

\vskip 2 pt

\centerline{$f(w\xi)=w^mf(\xi)$ and 
$f^{\sigma}(w \xi^{\sigma})=w^mf^{\sigma}(\xi^{\sigma}).$}

\vskip 2pt
\noindent
These spaces   are  respectively  invariant under 
$\pi_{\alpha}$ and   $\pi_{\alpha}^{\sigma}$. For  $f\in  {\cal O}_m(\Xi)$,  $h\in  {\cal O}_{m+{1\over 2}}(\Xi)$ and  $f^{\sigma}\in {\cal O}_m(\Xi^{\sigma})$, $h^{\sigma}\in {\cal O}_{m+{1\over 2}}(\Xi^{\sigma})$    their   corresponding functions
 $\phi$, $\psi$   and $\phi^{\sigma}$, $\psi^{\sigma}$  on ${\bboard C}^{r}$ satisfy for $\mu>0$,
 
 \vskip 2 pt
 
 \centerline{$\phi(\mu\cdot(z)=\mu^{{2m}}\phi(z)$,  \quad $\psi(\mu\cdot(z))=\mu^{{2m+1}}\psi(z)$,\quad }
 
\vskip 2 pt

\centerline{$\phi^{\sigma}(\mu\cdot(z))=\mu^{{2m}}\phi^{\sigma}(z)$, \quad $\psi^{\sigma}(\mu\cdot(z))=\mu^{{2m+1}}\psi^{\sigma}(z),$\quad }

 \vskip 2 pt
  \noindent
and  the correspondances 
$f(\xi)\mapsto \phi(z)$ and  $f^{\sigma}(\xi^{\sigma})\mapsto \phi^{\sigma}(z)$ 
map the spaces  ${\cal O}_m(\Xi)$, ${\cal O}_m(\Xi^{\sigma})$ and ${\cal O}_{m+{1\over 2}}(\Xi)$, ${\cal O}_{m+{1\over 2}}(\Xi^{\sigma})$,  to  respectively

\vskip 2 pt
 
 \centerline{${\cal O}_{2m}({\bboard C}^{r})=\{\phi\in {\cal O}({\bboard C}^{r}) \mid \phi(\mu z)=\mu^{{2m}}\phi(z)\},$}
 
 \vskip 4 pt
 
 \centerline{${\cal O}_{2m+1}({\bboard C}^{r})=\{\phi\in {\cal O}({\bboard C}^{r}) \mid \phi(\mu z)=\mu^{{2m+1}}\phi(z)\}$.}
 
\vskip 2 pt
 \noindent
 Let $\widetilde{\cal O}_{2m}({\bboard C}^{d_0})$,  $\widetilde{\cal O}_{2m+1}({\bboard C}^{d_0})$ and $\widetilde{\cal O}_{2m}^{\sigma}({\bboard C}^{d_0})$,   $\widetilde{\cal O}_{2m+1}^{\sigma}({\bboard C}^{d_0})$ be the sets of such functions $\phi$, $\phi^{\sigma}$  and $\psi$, $\psi^{\sigma}$  corresponding to the functions  $f\in{\cal O}_m(\Xi)$, $f^{\sigma}\in{\cal O}_m(\Xi^{\sigma})$, $h\in{\cal O}_{m+{1\over 2}}(\Xi)$, $h^{\sigma}\in {\cal O}_{m+{1\over 2}}(\Xi^{\sigma})$. 

\vskip 6 pt

Denote by  $\pi_{\alpha,m}$,  $\pi_{\alpha,m+{1\over 2}}$ and  $\pi_{\alpha,m}^{\sigma}$,  $\pi_{\alpha,m+{1\over 2}}^{\sigma}$
 the restrictions  of the representations $\pi_{\alpha}$ and $\pi_{\alpha}^{\sigma}$ to the spaces 
 ${{\cal O}}_m(\Xi)$, ${{\cal O}}_{m+{1\over 2}}(\Xi)$  and    ${{\cal O}}_m(\Xi^{\sigma})$, ${{\cal O}}_{m+{1\over 2}}(\Xi^{\sigma})$
  and  also  by $\pi_{\alpha,2m}$,   $ \pi_{\alpha,2m}^{\sigma}$,  $ \pi_{\alpha,2m+1}$, $ \pi_{\alpha,2m+1}^{\sigma}$ the corresponding representations on  $\widetilde{\cal O}_{2m}({\bboard C}^{r})$,  $\widetilde{\cal O}_{2m+1}({\bboard C}^{r})$, $\widetilde{\cal O}_{2m}^{\sigma}({\bboard C}^{r})$,   $\widetilde{\cal O}_{2m+1}^{\sigma}({\bboard C}^{r})$.
 It follows from Theorem 3.1. in [A.12b] that these spaces consist in polynomials and are finite dimensional,  ${{\cal O}}_{m}(\Xi)=\{0\}$ and ${{\cal O}}_{m}(\Xi^{\sigma})=\{0\}$ for $m<0$, and  
 $\pi_{\alpha,2m}$,   $\pi_{\alpha,2m+1}$,  $\pi_{\alpha,2m}^{\sigma}$, $\pi_{\alpha,2m+1}^{\sigma}$   are  irreducible.
\vskip 1 pt
\noindent
Also, 
  $L_{{\bboard R}}$-invariant norms on ${\cal O}_m(\Xi)$,   ${\cal O}_{m+{1\over 2}}(\Xi)$ and ${\cal O}_m(\Xi^{\sigma})$,   ${\cal O}_{m+{1\over 2}}(\Xi^{\sigma})$
  are  defined by:
\vskip 2pt
 \noindent

$$\Vert\phi\Vert_{m}^2={1\over a_{m}}\int_{{\bboard C}^{r}}\vert\phi(z)\vert^2H(z )^{-2m}m_0(dz),$$
$$\Vert\phi\Vert_{m+{1\over 2}}^2={1\over a_{m+{1\over 2}}}\int_{{\bboard C}^{r}}\vert\phi(z)\vert^2H(z)^{-2m+1}m_0(dz),$$

\noindent
where 
 \vskip 2 pt
 
 \centerline{$H(z)=\tau({1\over r}I_r+z\bar z)=1+{\rm tr}(z\bar z),$}

\vskip 2 pt
\noindent
and  $m_0(d(z))=H(z)^{-(r+1)}m(d(z)),$ is the $L$-invariant measure, 
 $p_m$  and $p_{m+{1\over 2}}$ are  suitable integers,   and   
$$a_{m}=\int_{{\bboard C}^{r}}H(z)^{-2m}m_0(dz), \quad a_{m+{1\over 2}}=\int_{{\bboard C}^{r}}H(z)^{-(2m+1)}m_0(dz). $$
\noindent
$$a_{m}=\pi^{r}{1\over (2m+d_0)\ldots(2m+1)}, \quad a_{m+{1\over 2}}=\pi^{r}{1\over (2m+1+d_0)\ldots(2m+2)}.$$
\vskip 2 pt
\noindent
These spaces 
 become invariant Hilbert spaces  with
reproducing kernels: 
$${\cal K}_{m}(\xi,\xi')=\Phi(\xi,\xi')^{2m} , \quad {\cal K}_{m+{1\over 2}}(\xi,\xi')
=\Phi(\xi,\xi')^{2m+1},$$
\noindent
with 

\vskip 2 pt

\centerline{$\Phi(\xi,\xi')=\tau({1\over r}I_r+z\bar z')$ for $\xi=\xi_{z}$, $\xi'=\xi_{z'}$.}

\vskip 2 pt
\noindent
They are   the irreducible $L_{{\bboard R}}$-invariant  subspaces of ${{\cal O}}(\Xi)$  and  ${{\cal O}}(\Xi^{\sigma})$.
The Fock spaces ${\cal F}(\Xi)$ and ${\cal F}(\Xi^{\sigma})$ are  Hilbert subspaces of ${\cal O}(\Xi)$ and ${\cal O}(\Xi)$ which are  $L$-invariant, they  therefore decompose
\vskip 2 pt

\centerline{${\cal F}(\Xi)=\sum\limits_{m=0}^{\infty}{\cal O}_m(\Xi)+\sum\limits_{m=0}^{\infty}{\cal O}_{m+{1\over 2}}(\Xi) ,$}

\vskip 2 pt

\centerline{${\cal F}(\Xi^{\sigma})=\sum\limits_{m=0}^{\infty}{\cal O}_m(\Xi^{\sigma})+\sum\limits_{m=0}^{\infty}{\cal O}_{m+{1\over 2}}(\Xi^{\sigma}) .$}

\vskip 4 pt
\noindent
The  Hilbert norms on ${\cal F}(\Xi)$ and ${\cal F}(\Xi)$  are   of the following form : for
\vskip 2 pt

\centerline{$f=\sum_{m=0}^{\infty}f_m+\sum_{m=0}^{\infty}f_{m+{1\over 2}} $}

\vskip 2 pt
\noindent

\vskip 2 pt
\noindent
then
\vskip 2 pt

\centerline{
$\Vert f\Vert_{\cal F}^2=\sum_{m=0}^{\infty}{1\over c_m}\Vert f_m\Vert_m^2+\sum_{m=0}^{\infty}{1\over c_{m+{1\over 2}}}\Vert f_{m+{1\over 2}}\Vert_{m+{1\over 2}}^2,$}

\vskip 2 pt
\noindent

\noindent
The sequences $(c_m)$, $(c_{m+{1\over 2}})$,  are determined in such a way that the representation $\rho_{\alpha}$ of ${\goth g}_{\bboard R}$ and $\rho_{\alpha}^{\sigma}$ are  unitary (see [A12b], Theorem 5.1).   One gets
\vskip 2 pt

\centerline{$c_m={1\over(2m)!}$, \quad $c_{m+{1\over 2}}={1\over(2m+1)!},$}

\vskip 2 pt
\noindent
and then  ${\cal F}(\Xi)$ and ${\cal F}(\Xi)$ turn out to be the classical Fock space ${\cal F}({\bboard C}^r).$
\vskip 2 pt
\hfill
\eject
\noindent
For the representations $\rho_{\alpha}$ and $\rho_{\alpha}^{\sigma}$  of ${\goth g}$,  the elements  $\omega \in {\goth l}$  act by  
$$\rho_{\alpha}(\omega)={\rm d}\pi_{\alpha}(\omega)-{1\over 2}d\pi_{\alpha}(H_0)$$
\noindent
 and 
 $$\rho_{\alpha}^{\sigma}(\omega)={\rm d}\pi_{\alpha}^{\sigma}(\omega)-{1\over 2}d\pi_{\alpha}^{\sigma}(H_0).$$

\vskip 4 pt
  \noindent
  For  $\rho_{\alpha}$, the  elements  $p\in{\cal V}$  act  by multiplication and  the elements $p^{\sigma}\in{\cal V}^{\sigma}$   act by differentiation, and, 
for $\rho_{\alpha}^{\sigma}=\pi(\sigma)\rho_{\alpha}\pi(\sigma)$, the elements  $p\in{\cal V}$  act  by differentiation  and    the  elements  $p^{\sigma}\in{\cal V}^{\sigma}$  act by multiplication.
\vskip 2 pt
\noindent
The representation $\rho_{\alpha}$ is determined by the operators $\rho_{\alpha}(E)$ which involves the multiplication  operator $\tau(z^2)$, and $\rho_{\alpha}(F)=-\rho_{\alpha}(E)^*$.  The representation $\rho_{\alpha}^{\sigma}$ is determined by the operators $\rho_{\alpha}^{\sigma}(E)$ which involves the differential   operator $\tau({\partial^2\over\partial z^2})$, and $\rho_{\alpha}(F)=-\rho_{\alpha}(E)^*$. 
\vskip 2 pt
\noindent
The operators  $\rho_{\alpha}(E), \rho_{\alpha}(F), \rho_{\alpha}(H_0)$ and $\rho_{\alpha}^{\sigma}(E), \rho_{\alpha}^Ô{\sigma}(F), \rho_{\alpha}^{\sigma}(H_0)$
 are     given by
$$\rho_{\alpha}(E)\phi(z)={i\over 4}\tau(z^2)\phi(z),$$
$$\rho_{\alpha}^{\sigma}(E)\phi(z)={i\over 4}\tau({\partial^2\over\partial z^2})\phi(z) $$
$$\rho_{\alpha}(F)\phi(z)={i\over 4}\tau({\partial^2\over\partial z^2})\phi(z),$$
$$\rho_{\alpha}^{\sigma}(F)\phi(z)={i\over 4}\tau(z^2)\phi(z),$$
$$\rho_{\alpha}({H_0})\phi(z,z')=(1-r)\bigr(-\alpha r\phi(z,z')+{1\over 2}{\cal E}\phi(z,z')\bigl),$$
$$\rho_{\alpha}^{\sigma}({H_0})\phi(z,z')=(1-r)\bigr(\alpha r\phi(z,z')+{1\over 2}{\cal E}\phi(z,z')\bigl),$$

\vskip 6 pt
\noindent
 where 
$\alpha=-{1\over 4}$ 
  and  ${\cal E}$ is the Euler operator given by 
$$({\cal E}\phi)(z)={{\rm d}\over {\rm d}s}_{\mid_{s=1}}\phi(sz).$$
\vskip 4 pt
\hfill
\eject
\noindent
\bf 2. Some harmonic analysis\rm  
\vskip 6 pt
\noindent
We consider on $V={\rm Sym}(r,{\bboard C})$ the homogeneous polynomial
\vskip 4 pt

\centerline{$Q(x)=\det(x)^2.$}

\vskip 2 pt
\noindent
Then the structure group is 
$$L={\rm Str}(V,\Delta)={\rm GL}(r,{\bboard C})$$
\noindent
(quotiented by $\{\pm I_r)\}$). 
The orbit $\Xi$ of $\tau_{\lambda}$ under $L$  has  dimension $r$ and  can be identified to 
\vskip 4 pt

\centerline{$\Gamma(\Xi)=\{\xi_z=zz^t \mid z\in {\bboard C}^r\}$.}

\vskip 4 pt
\noindent
Let ${\cal Y}_m({\bboard R}^{r})$ be the space of spherical harmonics of degree $m$ on ${\bboard R}^{r}$: harmonic polynomials which are homogeneous of degree $m$ on ${\bboard R}^{r}$. The map
$${\cal Y}_m({\bboard R}^{r}) \rightarrow {\cal O}_m(\Xi), \quad \Phi \mapsto f,$$
\noindent
given by 
$$f(\xi_z)=\int_S<z,x>^m\Phi(x)s(dx),$$
\noindent
where  
$$<z,x>=\sum\limits_{j=1}^{r}z_jx_j,$$
\noindent
is an isomorphism which intertwins the representations of $O(r)$  on both spaces. ($S$ is the unit sphere in ${\bboard R}^{r}$, and $s(dx)$ is the uniform  measure on $S$ with total measure equal to one (see [F.15], section 2).

\bigskip
For a holomorphic function $a$ on ${\bboard C}$ we define the integral operator  $A$ from  ${\cal C}(S)$ into ${\cal O}(\Xi)$:
$$Af(\xi_z)=\int_Sa(<z,x>)f(x)s(dx).$$
\noindent
The operator $A$ is equivariant  with respect to the action of $O(r)$ and maps ${\cal Y}_m({\bboard R}^{r})$ into ${\cal O}_m(\Xi)$. 

\bigskip
\hfill
\eject
\noindent
\bf 3. The Lie algebra ${\goth g}$ and its isomorphism with ${\goth sp}(r,{\bboard C})$\rm 
\vskip 6 pt
\noindent
The Lie algebra ${\goth g}$ is isomorphic to  the matrix Lie algebra 
 \vskip 1 pt

 \centerline{$\widetilde{\goth g}=\Big\{\pmatrix{\omega_1&v\cr
 u&-\omega_1^t} \mid u, v\in {\rm Sym}(r,{\bboard C}), \quad \omega_1\in{\goth gl}(r,{\bboard C})\Big\}$, }
 
 \noindent
and the isomorphism is given by 
$$\tau_u+\omega+\tau_v^{\sigma}\in {\goth g}={\cal V}\oplus{\goth l}\oplus{\cal V}^{\sigma} \mapsto \widetilde{\tau_u}+\widetilde{\omega}+\widetilde{\tau_v^{\sigma}}\in \widetilde{\goth g}$$
\noindent
with $\omega=(\omega_1,-\omega_1^t)$, 
  \vskip 2 pt
  
  \centerline{$ \widetilde{\omega}=\pmatrix{\omega_1&0\cr
 0&-\omega_1^t}+{\rm tr}(\omega_1)\pmatrix{-I_r&0\cr
 0&I_r},$}
 
 \vskip 2 pt

  \centerline{$\widetilde{\tau_u}={1\over 2}\pmatrix{0&0\cr
 u&0},  \quad  \widetilde{\tau_v^{\sigma}}={1\over 2}\pmatrix{0&v^t\cr
 0&0}$.}
 
 \vskip 4 pt
 \noindent
Then, ${\goth g}$ is isomorphic to ${\goth sp}(r,{\bboard C})$.
In fact this follows from  Proposition 1.1. in[A12b]:  every  $\omega=(\omega_1,-\omega_1^t)\in {\goth l}$  acts  on $V$ by $\omega x=\omega_1\cdot x+x\cdot\omega_1^t$. It follows that  for every $\omega, \omega'\in \goth l$, one has
$[\widetilde\omega,\widetilde{\omega'}]=\widetilde{[\omega,\omega']}.$
Moreover, for $\omega=(\omega_1,-\omega_1^t) \in {\goth l}$, one has
\vskip 4 pt

\centerline{$[\omega,\tau_u]=d\kappa(\omega)\tau_u=\tau_{-\omega_1^tu-u\omega_1+2({\rm tr}(\omega_1))u}$}

\vskip 1 pt
\noindent
and
\vskip 1 pt

\centerline{$[\widetilde \omega,\widetilde{\tau_u}]={1\over 2}\pmatrix{0&0\cr
-\omega_1^tu-u\omega_1+2{\rm tr}(\omega_1)u&0}$}

\vskip 1 pt
\noindent
then 
\vskip 1 pt

\centerline{$[\widetilde \omega,\widetilde{\tau_u}]=\widetilde{[\omega,\tau_u]}.$}

\vskip 4 pt
\noindent
One has also 
\vskip 1 pt

\centerline{$[\omega,\tau_v^{\sigma}]=d\kappa(\omega)\tau_v^{\sigma}=\tau_{\omega_1^tv+v\omega_1^t+2{\rm tr}(\omega_1)v}$}

\vskip 1 pt
\noindent
and

\centerline{$[\widetilde \omega,\widetilde{\tau_v^{\sigma}}]={1\over 2}\pmatrix{0&v^t\omega_1^t+\omega_1v^t+2{\rm tr}(\omega_1)v^t\cr
0&0}$}

\vskip 1 pt
\noindent
then 
\vskip 1 pt

\centerline{$[\widetilde \omega,\widetilde{\tau_v^{\sigma}}]=\widetilde{[\omega,\tau_v^{\sigma}]}.$}

\vskip 4 pt
\noindent
The matrices corresponding to $E=\tau$, $F=\tau^{\sigma}$ and $H_0$ (in (i) and (ii))  are 
\vskip 1 pt

\centerline{$\widetilde E={1\over 2}\pmatrix{0&0\cr
I_r&0}$, \quad $\widetilde F={1\over 2}\pmatrix{0&I_r\cr
0&0},  \widetilde{H_0}=(1-r)\pmatrix{-{1\over 2}I_r&0\cr
0&{1\over 2}I_r}$}

\noindent
then $[\widetilde E,\widetilde F]={1\over 2(1-r)}\widetilde{H_0}$ and,  since $[E,F]=(\lambda_0(E,F)+{1\over 2})H_0=({1\over 2}-{r\over 2(r-1)})H_0$, then $[\widetilde E,\widetilde F]=\widetilde{[E,F]}.$
Observe that  

\vskip 1 pt

\centerline{$[\widetilde H_0,\widetilde E]=-(r-1)\widetilde E$ and $[\widetilde H_0,\widetilde F]=(r-1)\widetilde F$, }

\vskip 1 pt
\noindent
i.e. 
 \vskip 1 pt
 
 \centerline{$[\widetilde H_0,\widetilde E]=\widetilde{[H_0,E]}$ and $[\widetilde H_0,\widetilde F]=\widetilde{[H_0,F]}$.}
 
\vskip 3 pt
\noindent
This proves that we have obtained an explicit Lie algebra isomorphism from ${\goth g}$ to  $\widetilde{\goth g}={\goth sp}(r,{\bboard C})$. 
\vskip 6 pt
\noindent
Now, let's consider the image $\widetilde{{\goth g}_{\bboard R}}$ of the real form ${\goth g}_{\bboard R}$ by this isomorphism. Since  ${\goth g}_{\bboard R}$ is given by  ${\goth g}_{\bboard R}={\goth l}_{\bboard R}+{\cal W}_{\bboard R}$, where ${\goth l}_{\bboard R}$ is the compact real form of  ${\goth l}$ and ${\cal W}_{\bboard R}$ is generated by the elements $\tau_u+\tau_u^{\sigma}$ and $i(\tau_u-\tau_u^{\sigma})$ for $u\in V_{\bboard R}={\rm Sym}({\bboard R})$, then  $\widetilde{{\goth g}_{\bboard R}}$ is  given by $\widetilde{{\goth g}_{\bboard R}}=\widetilde{{\goth l}_{\bboard R}}+\widetilde{{\cal W}_{\bboard R}}$ with
$$\eqalignno{\widetilde{{\goth l}_{\bboard R}}=\{\pmatrix{\omega_1&0\cr
0&-\omega_1^t} \mid  \omega_1=w_1+iw_1',   w_1\in {\rm Skew}(r,{\bboard R}),  w_1' \in {\rm Sym}(r,{\bboard R })\},}$$
$$\widetilde{{\cal W}_{\bboard R}}=\{\pmatrix{0&u-iv\cr
u+iv&0} \mid u,v\in {\rm Sym}(r,{\bboard R})\},$$
\noindent
i.e.
$$\widetilde{{\goth g}_{\bboard R}}=\{\pmatrix{\omega_1&u-iv\cr
u+iv&-\omega_1^t} \mid \omega_1=w_1+iw_1'\in {\goth su}(r),   u,v  \in {\rm Sym}(r,{\bboard R })\}.$$
\vskip 2 pt
\noindent 
Let $g_0$ be the element of  ${\rm Sp}(r,{\bboard C})$ given by
$$g_0={1\over\sqrt 2}\pmatrix{iI_r&I_r\cr
-I_r&-iI_r}.$$
\noindent
Since 
$${\rm Ad}(g_0)\pmatrix{w_1&w_1'\cr
-w_1'&w_1}=\pmatrix{w_1+iw_1'&0\cr
0&w_1-iw_1'}$$ 
and 
$${\rm Ad}(g_0)\pmatrix{v&-u\cr
-u&-v}=\pmatrix{0&u-iv\cr
u+iv&0},$$
\noindent
then 
\vskip 4 pt
\noindent

\noindent

\centerline{${\rm Ad}(g_0^{-1})(\widetilde{\goth l}_{\bboard R})=\{\pmatrix{w_1&w_1'\cr
-w_1'&w_1} \mid    w_1\in {\rm Skew}(r,{\bboard R}),  w_1' \in {\rm Sym}(r,{\bboard R })\},$}
 
 \vskip 2 pt
 \noindent
 and
 \vskip 4 pt
 
 \centerline{${\rm Ad}(g_0^{-1})(\widetilde{\cal W}_{\bboard R})=\{\pmatrix{v&-u\cr
-u&-v} \mid u,v\in {\rm Sym}(r,{\bboard R})\},$}

\vskip 8 pt
\noindent
in such a way that
$$\eqalignno{&{\rm Ad}(g_0^{-1})(\widetilde{{\goth g}_{\bboard R}})=\cr
&\{\pmatrix{w_1+v&w_1'-u\cr
-w_1'-u&w_1-v},   w_1\in {\rm Skew}(r,{\bboard R}),  w_1', u,v \in {\rm Sym}(r,{\bboard R })\}
={\goth sp}(r,{\bboard R}).}$$

\hfill
\eject
Furthertmore, it  is well known that the minimal nilpotent $\widetilde G$-orbit in $\widetilde{\goth g}$ is 
\vskip 2 pt

\centerline{$\widetilde{\cal O}_{\rm min}={\rm Ad}(\widetilde G)(e_{11})$}

\vskip 2 pt
\noindent
where  $E_{11}$ is the diagonal matrix $E_{11}={\rm diag}(1,0,\ldots,0)$. and $e_{11}=\pmatrix{0&E_{11}\cr
0&0}=2\widetilde{\tau_{c_{\lambda}}^{\sigma}}$  is the highest root. It follows that

\vskip 4 pt

\centerline{$\widetilde{{\cal O}_{\rm min}}\cap \widetilde{{\cal W}}=\widetilde{{\cal O}_{\rm min}}\cap \widetilde{{\cal V}}\cup \widetilde{{\cal O}_{\rm min}}\cap \widetilde{{\cal V}^{\sigma}}$}

\vskip 4 pt
\noindent
with
$$\eqalign{\widetilde{{\cal O}_{\rm min}}\cap \widetilde{{\cal V}^{\sigma}}&={\rm Ad}(\widetilde L)(e_{11})\cr
&=\{{\rm Ad}(l)(e_{11}) \mid  l=\pmatrix{l_1&0\cr
0&(l_1^t)^{-1}} \mid l_1\in  {\rm GL}(r,{\bboard C})\}\cr
&=\{\pmatrix{0&l_1E_{11}l_1^t\cr
0&0} \mid l_1\in  {\rm GL}(r,{\bboard C})\}\cr
&=\{\pmatrix{0&zz^t\cr
0&0} \mid z\in {\bboard C}^r\}},$$
\noindent
and
$$\eqalign{\widetilde{{\cal O}_{\rm min}}\cap{\cal V}&={\rm Ad}(\widetilde L)({\rm Ad}(J)e_{11})\cr
&\{{\rm Ad}(l)({\rm Ad}(J)e_{11}) \mid l=\pmatrix{l_1&0\cr
0&(l_1^t)^{-1}} \mid l_1\in  {\rm GL}(r,{\bboard C})\}\cr
&=\{\pmatrix{0&0\cr
l_1E_{11}l_1^t&0} \mid l_1\in  {\rm GL}(r,{\bboard C})\}\cr
&=\{\pmatrix{0&0\cr
zz^t&0} \mid z\in {\bboard C}^r\}.}$$
\noindent
It follows that $\widetilde{{\cal O}_{\rm min}}\cap \widetilde{{\cal V}^{\sigma}}$ and $\widetilde{{\cal O}_{\rm min}}\cap{\cal V}$ are respectively the images of the orbits $\Xi$ and $\Xi^{\sigma}$ by the isomorphism ${\goth g}\rightarrow \widetilde{\goth g}$ and are  diffeomorphic  to 
$$\Gamma(\Xi)=\{\xi_z=zz^t \mid z\in {\bboard C}^r\}$$
\vskip 1 pt
\noindent
and that the  map $\kappa(\sigma)$ corresponds here to  $X\mapsto {\rm Ad}(J)X$.
\vskip 4 pt
\vskip 2 pt
\noindent
Moreover, 

\vskip 4 pt

\centerline{$\widetilde{{\cal O}_{\rm min}}\cap {\goth sp}(r,{\bboard R})=\widetilde{{\cal O}_{\rm min}}\cap{\rm Ad}(g_0^{-1})(\widetilde{{\goth g}_{\bboard R}})=Y^+\cup Y^-.$}

\vskip 4 pt
\noindent
\vskip 4 pt
\noindent
\vskip 2 pt
\noindent

\hfill
\eject
\bigskip
\noindent
\bf 4. The Schr\" odinger  model. \rm 

\bigskip

\noindent
Let $\Gamma_{\bboard R}$ be the  open  cone in ${\bboard R}^{r}$ and  $S$ be the unit sphere  given respectively by:
$$\Gamma_{\bboard R}=\{x\in {\bboard R}^{r} \mid \vert x\vert\ne 0\},$$
\noindent
and
$$S=\{x\in {\bboard R}^{r} : \vert x\vert=1\}.$$
\vskip 6 pt
\noindent
The group $L_{\bboard R}={\rm GL}(r,{\bboard R})$ acts on ${\bboard R}^{r}$ by the natural  representation.
 This action stabilizes the cone $\Gamma_{\bboard R}$. The multiplicative group ${\bboard R}_+^*$ acts on $\Gamma_{\bboard R}$ as a dilation and the quotient space $M=\Gamma_{\bboard R}/{\bboard R}_+^*$ is identified with $S$. This defines  an action of $L_{\bboard R}$ on $S$,  which leads to a $L_{\bboard R}$-equivariant principal ${\bboard R}_+^*$-bundle:
$ \Gamma_{\bboard R} \rightarrow S, x\mapsto {x\over\vert x\vert}.$ 
For $\lambda \in {\bboard C}$, let 
${\cal E}_{\lambda}(\Gamma_{\bboard R})$ and  ${\cal E}_{\lambda}^{\sigma}(\Gamma_{\bboard R})$  be  the spaces of ${\cal C}^{\infty}$-functions on $\Gamma_{\bboard R}$ homogeneous  of degree $\lambda$:
$${\cal E}_{\lambda}(\Gamma_{\bboard R})={\cal E}_{\lambda}^{\sigma}(\Gamma_{\bboard R})=\{f\in {\cal C}^{\infty}(\Gamma_{\bboard R}) \mid f(tx)=t^{\lambda}f(x), \quad x\in \Gamma, t >0\},$$
\noindent
The group $L_{\bboard R}={\rm GL}(r,{\bboard R})$ acts naturally on ${\cal E}_{\lambda}(\Gamma_{\bboard R})$, and, under  the action of the subgroup $O(r)$, the space   ${\cal E}_{\lambda}(\Gamma_{\bboard R})$ decomposes as:
$${\cal E}_{\lambda}(\Gamma_{\bboard R})\mid_S \simeq \bigoplus_{k=0}^{\infty}{\cal Y}_{k}({\bboard R}^{r}).$$

\vskip 4 pt
\noindent
These representations  extend to representations $R$ and $R^{\sigma}$ of the metaplectic  group ${\rm Mp}(r,{\bboard R})$ on  the Hilbert space ${\rm L}^2({\bboard R}^r)$  as follows: 
denote by 
\vskip 6 pt

\centerline{$g(l_1)=\pmatrix{l_1&0\cr
0&(l_1^t)^{-1}} $ for $l_1\in {\rm GL}(r,{\bboard R})$,}

\vskip 6 pt

\centerline{$t(u)=\exp(2\widetilde{\tau_u})=\pmatrix{1&0\cr
u&1}$, $t(u)^{\sigma}=\exp(2\widetilde{\tau_u^{\sigma}})=\pmatrix{1&u\cr
0&1}$, ($u\in {\rm Sym}(r)$),}

\vskip 6 pt

\centerline{$J_r=\pmatrix{0&-I_r\cr
I_r&0}$.}

\vskip 6 pt
\noindent
It is well  known that the elements $g(l_1), t(u)$ and $J_r$, for $l_1\in {\rm GL}(r,{\bboard R}), u\in  {\rm Sym}(r,{\bboard R})$, generate  the group symplectic group  ${\rm Sp}(r,{\bboard R})$. One considers the  two  representations $R$ and $R^{\sigma}$  of the metaplectic group ${\rm Mp}(r,{\bboard R})$ on  ${\rm L}^2({\bboard R}^r)$,   determined by:
\hfill
\eject
$$R(g(l_1))f(x)=(\det l_1)^{{1\over 2}}f(l_1^tx),$$
$$R(t(u)^{\sigma})f(x)=e^{-{i\over 2}\tau_u(x^2)}f(x),$$
$$R(J_r)f(x)=a_0\int_{{\bboard R}^r}e^{i\tau(xy)}f(y)dy.$$
\noindent
and
$$R^{\sigma}(g(l_1))f(x)=(\det l_1)^{-{1\over 2}}f(l_1^tx),$$
$$R^{\sigma}(t(u))f(x)=e^{-{i\over 2}\tau_u(x^2)}f(x),$$
$$R^{\sigma}(J_r)f(x)=a_0\int_{{\bboard R}^r}e^{i\tau(xy)}f(y)dy.$$
\noindent
Then
$$dR(\widetilde{\tau_u^{\sigma}})f(x)=-{i\over 4}\tau_u(x^2)f(x),\quad dR(\widetilde{\tau_u})f(x)=-{i\over 4}\tau_u({\partial^2\over\partial x^2})f(x)$$
\noindent
and 
$$dR^{\sigma}(\widetilde{\tau_u})f(x)=-{i\over 4}\tau_u(x^2)f(x),\quad dR^{\sigma}(\widetilde{\tau_u^{\sigma}})f(x)=-{i\over 4}\tau_u({\partial^2\over\partial x^2})f(x)$$
\noindent
In particular
\vskip 2 pt
$$dR(\widetilde E)f(x)=-{i\over 4}\tau({\partial^2\over\partial x^2})f(x),\quad dR(\widetilde F)f(x)=-{i\over 4}\tau(x^2)f(x)$$
\noindent
and
$$dR^{\sigma}(\widetilde F)f(x)=-{i\over 4}\tau({\partial^2\over\partial x^2})f(x),\quad dR^{\sigma}(\widetilde E)f(x)=-{i\over 4}\tau(x^2)f(x)$$
\vskip 4 pt
\noindent
Denote by 
$${\rm L}^2({\bboard R}^r)_{\rm even}=\{f\in {\rm L}^2({\bboard R}^r) \mid f(-x)=f(x)\}$$
\noindent
and
$${\rm L}^2({\bboard R}^r)_{\rm odd}=\{f\in {\rm L}^2({\bboard R}^r)\mid f(-x)=-f(x)\}.$$

\vskip 4 pt
\noindent
The following facts are well-known:
\vskip 4 pt
\noindent
1) (Irreducibility) The representations   $\bigr(R,{\rm L}^2({\bboard R}^r)_{\rm even}\bigl)$,    $\bigr(R,{\rm L}^2({\bboard R}^r)_{\rm odd}\bigl)$, $\bigr(R^{\sigma},{\rm L}^2({\bboard R}^r)_{\rm even}\bigl)$,    $\bigr(R^{\sigma},{\rm L}^2({\bboard R}^r)_{\rm odd}\bigl)$ of  ${\rm Mp}(r,{\bboard R})$  are irreducible.

\vskip 2 pt
\noindent
2) ($K$-type decomposition) The  underlying  $({\goth g},K)$-modules,  $(R)_{K^L}$ and $(R_{\alpha_0}^{\sigma})_{K^L}$, for $K=K^L=O(r,{\bboard C})$,  have the following  $K$-type formulas 
\vskip 4 pt
\centerline{
$(R)_{K^L}=\bigoplus\limits_{m=0}^{\infty}{\cal Y}_{2m}({\bboard R}^{r})+\bigoplus\limits_{m=0}^{\infty}{\cal Y}_{2m+1}({\bboard R}^{r})$,}

\vskip 2 pt

\centerline{ 
$(R^{\sigma})_{K^L}=\bigoplus\limits_{m=0}^{\infty}{\cal Y}_{2m}({\bboard R}^{r})+\bigoplus\limits_{m=0}^{\infty}{\cal Y}_{2m+1}({\bboard R}^{r}).$}

\vskip 2 pt
\noindent
3) (Unitarity) The representations $R$ and $R^{\sigma}$  of ${\rm Mp}(r,{\bboard R})$ on ${\rm L}^2({\bboard R}^r)$
 are unitary.
\hfill
\eject
\noindent
\bf 5. The  intertwining operator\rm 
\vskip 4 pt
\noindent
Recall from section 3 that
$$\eqalignno{&{\rm Ad}(g_0^{-1})(\widetilde{{\goth g}_{\bboard R}})=\cr
&\{\pmatrix{w_1+v&w_1'-u\cr
-w_1'-u&w_1-v},   w_1\in {\rm Skew}(r,{\bboard R}),  w_1', u,v \in {\rm Sym}(r,{\bboard R })\}
={\goth sp}(r,{\bboard R})}$$
\noindent
where 
$$g_0={1\over\sqrt 2}\pmatrix{iI_r&I_r\cr
-I_r&-iI_r}.$$
\vskip 2 pt
\noindent
Denote by $\widetilde{G}={\rm Sp}(r,{\bboard C})$.  It is generated by the elements $g(l_1), t(u)$ and $J_r$ (for $l_1\in {\rm GL}(r,{\bboard C}), u\in  {\rm Sym}(r,{\bboard C})$).
\vskip 2 pt
\noindent
The representations $\rho:=\rho_{\alpha}$ and  $\rho^{\sigma}:=\rho_{\alpha}^{\sigma}$  in section 1(i.e with $\alpha=-{1\over 4}$)  'integrate'  to  the representations $T$  and $T^{\sigma}$ of  $\widetilde{G}={\rm Mp}(r,{\bboard C})$ on ${\cal O}_{\rm fin}(\Xi)$  given by
$$T(g(l_1))\phi(z)=(\det l_1)^{2\alpha}\phi(l_1^{-1}z),$$
$$T(t(u))\phi(z)=e^{{i\over 2}\tau_u(z^2)}\phi(z),$$
$$T(J)\phi(z)=a_0\int_{{\bboard R}^r}e^{i\tau(zy)}\phi(y)dy$$
\noindent
and
$$T^{\sigma}(g(l_1))\phi(z)=(\det l_1)^{-2\alpha}\phi(l_1^tz),$$
$$T^{\sigma}(t(u)^{\sigma})\phi(z)=e^{{i\over 2}\tau_u(z^2)}\phi(z),$$
$$T^{\sigma}(J)\phi(z)=a_0\int_{{\bboard R}^r}e^{i\tau(zy)}\phi(y)dy$$
\noindent
in such a way that
$$dT(\widetilde{\tau_u^{\sigma}})\phi(z)={i\over 4}\tau_u(z^2)\phi(z),\quad dT(\widetilde{\tau_u})\phi(z)={i\over 4}\tau_u({\partial^2\over\partial z^2})\phi(z),$$
$$dT^{\sigma}(\widetilde{\tau_u})\phi(z)={i\over 4}\tau_u(z^2)\phi(z),\quad dT^{\sigma}(\widetilde{\tau_u^{\sigma}})\phi(z)={i\over 4}\tau_u({\partial^2\over\partial z^2})\phi(z)$$
\noindent
(where we precise that the exponant $2\alpha$ arises  from $\chi(l)^{\alpha}=(\det l_1)^{2\alpha}$).
\vskip 2 pt
\hfill
\eject
\noindent
Denote  by 

\vskip 4 pt

\centerline{$\widetilde{G_{\bboard R}}={\rm Ad}(g_0)({\rm Mp}(r,{\bboard R})$.}

\vskip 6 pt
\noindent
and consider the  Hilbert spaces $T_0^{-1}({\cal F}({\bboard C}^r))$  and ${T_0^{\sigma}}^{-1}({\cal F}({\bboard C}^r))$, equipped with the norms 

\vskip 6 pt

\centerline{$\Vert\psi\Vert_{T_0^{-1}({\cal F}({\bboard C}^r))}=\Vert T_0\psi\Vert_{{\cal F}({\bboard C}^r)}$}

\vskip 6 pt
\noindent
and
\vskip 6 pt

\centerline{$\Vert\psi\Vert_{{T_0^{\sigma}}^{-1}({\cal F}({\bboard C}^r))}=\Vert T_0^{\sigma}\psi\Vert_{{\cal F}({\bboard C}^r)}.$}

\vskip 6 pt
\noindent
Then, $\bigr(T,T(g_0^{-1})({\cal F}({\bboard C}^r))\bigl)$  and  $\bigr(T^{\sigma},T^{\sigma}(g_0^{-1})({\cal F}({\bboard C}^r))\bigl)$
 are unitary representations of  the metaplectic group ${\rm Mp}(r,{\bboard R})$. In fact,  for $g\in {\rm Mp}(r,{\bboard R})$, there is $g'\in \widetilde{G_{\bboard R}}$ such that
 
 \vskip 4 pt
 
 \centerline{$g=g_0^{-1}g'g_0$, }
 
 \vskip 4 pt
 \noindent
 then
 \vskip 4 pt
 
 \centerline{$T(g)=T(g_0^{-1}T(g')T(g_0)$ and  $T^{\sigma}(g)=T^{\sigma}(g_0^{-1}T^{\sigma}(g')T^{\sigma}(g_0)$.}
 
 \vskip 4 pt
 \noindent
 It follows that for  $\psi\in T_0^{-1}({\cal F}({\bboard C}^r))$,   $\psi^{\sigma}\in (T_0^{\sigma})^{-1}({\cal F}({\bboard C}^r))$ and for $g\in {\rm Mp}(r,{\bboard R})$,
$$\eqalign{&\Vert T(g)\psi\Vert_{T(g_0^{-1})({\cal F}({\bboard C}^r))}=\Vert T(g_0^{-1})T(g')T(g_0)\psi\Vert_{T(g_0^{-1})({\cal F}({\bboard C}^r))}\cr
&=\Vert T(g')T(g_0)\psi\Vert_{{\cal F}({\bboard C}^r)}=\Vert T(g_0)\psi\Vert _{{\cal F}({\bboard C}^r)}=\Vert\psi\Vert_{T(g_0^{-1})({\cal F}({\bboard C}^r))}}$$
\noindent
and similarily,
 $$\eqalign{&\Vert T^{\sigma}(g)\psi^{\sigma}\Vert_{T^{\sigma}(g_0^{-1})({\cal F}({\bboard C}^r))}=\Vert T^{\sigma}(g_0^{-1})T^{\sigma}(g')T^{\sigma}(g_0)\psi^{\sigma}\Vert_{T^{\sigma}(g_0^{-1})({\cal F}({\bboard C}^r))}\cr
&=\Vert T^{\sigma}(g')T^{\sigma}(g_0)\psi^{\sigma}\Vert_{{\cal F}({\bboard C}^r)}=\Vert T^{\sigma}(g_0)\psi^{\sigma}\Vert _{{\cal F}({\bboard C}^r)}=\Vert\psi^{\sigma}\Vert_{T^{\sigma}(g_0^{-1})({\cal F}({\bboard C}^r))}.}$$
\vskip 4 pt
\noindent
Recall that the Bargmann transform 
$${\cal B} : L^2({\bboard R}^r) \rightarrow {\cal F}({\bboard C}^r)$$
\noindent
is a unitary operator given by the integral formula
$$({\cal B}f)(z)=\pi^{-{r\over 4}}\int_{{\bboard R}^r}e^{-{1\over 2}(\tau(x^2)+\tau(z^2))+\sqrt 2\tau(zx)}f(x)dx$$
\noindent
where $dx$ is Lebesgue measure, and $\tau$ is the usual  symmetric bilinear form in $r$ variables.
\vskip 2 pt
\hfill
\eject
\noindent
\bf Theorem   5.1| \rm  
\vskip 2 pt
\noindent
(i) The unitary  representations $\bigr(R,L^2({\bboard R}^r)\bigl)$  and $\bigr(T^{\sigma},T^{\sigma}(g_0^{-1})({\cal F}({\bboard C}^r)\bigl)$ of the group ${\rm Mp}(r,{\bboard R})$ are unitarily equivalent. The intertwinnig operator is given by ${\cal B}_0=T(g_0^{-1})\circ{\cal B}$. 
\vskip 2 pt
\noindent
(ii) The unitary  representations $\bigr(R^{\sigma},L^2({\bboard R}^r)\bigl)$  and $\bigr(T,T(g_0^{-1})({\cal F}({\bboard C}^r)\bigl)$ of the group ${\rm Mp}(r,{\bboard R})$ are unitarily equivalent. The intertwinnig operator is given by ${\cal B}_0^{\sigma}=T^{\sigma}(g_0^{-1})\circ{\cal B}$. 

\vskip 2 pt
\noindent
\bf Proof. \rm 
\vskip 2 pt
\noindent
In fact, if an operator ${\cal B}_0: L^2({\bboard R}^r) \rightarrow T^{\sigma}(g_0^{-1})({\cal F}({\bboard C}^r)$ intertwins $R$ and $T^{\sigma}$, then for every $g=g_0^{-1}g'g_0\in {\rm Mp}(r,{\bboard R})$, one has
\vskip 6 pt

\centerline{$T^{\sigma}(g){\cal B}_0={\cal B}_0R(g),$}

\vskip 4 pt
\noindent
i.e.
\vskip 4 pt

\centerline{$T^{\sigma}(g_0^{-1})T^{\sigma}(g')T^{\sigma}(g_0){\cal B}_0={\cal B}_0R(g)$}

\vskip 4 pt
\noindent
which means
\vskip 4 pt

\centerline{$T^{\sigma}(g')(T^{\sigma}(g_0){\cal B}_0)=(T^{\sigma}(g_0){\cal B}_0)R(g)$}

\vskip 4 pt
\noindent
i.e.
\vskip 4 pt

\centerline{$T^{\sigma}(g')\widetilde{\cal B}=\widetilde{\cal B}R(g)$ with $\widetilde{\cal B}=(T^{\sigma}(g_0){\cal B}_0)$.}

\vskip 6 pt
\noindent
Similarily,  if an operator ${\cal B}_0^{\sigma}: L^2({\bboard R}^r) \rightarrow T(g_0^{-1})({\cal F}({\bboard C}^r)$ intertwins $R^{\sigma}$ and $T$, then for every $g=g_0^{-1}g'g_0\in {\rm Mp}(r,{\bboard R})$, one has
\vskip 6 pt

\centerline{$T(g){\cal B}_0^{\sigma}={\cal B}_0^{\sigma}R^{\sigma}(g),$}

\vskip 4 pt
\noindent
i.e.
\vskip 4 pt

\centerline{$T(g_0^{-1})T(g')T(g_0){\cal B}_0^{\sigma}={\cal B}_0^{\sigma}R^{\sigma}(g)$}

\vskip 4 pt
\noindent
which means
\vskip 4 pt

\centerline{$T(g')(T(g_0){\cal B}_0^{\sigma})=(T(g_0){\cal B}_0^{\sigma})R^{\sigma}(g)$}

\vskip 4 pt
\noindent
i.e.
\vskip 4 pt

\centerline{$T(g')\widetilde{\cal B}^{\sigma}=\widetilde{\cal B}^{\sigma}R^{\sigma}(g)$ with $\widetilde{\cal B}^{\sigma}=(T(g_0){\cal B}_0^{\sigma})$.}

\bigskip

In what follows, we will see that $\widetilde{\cal B}=\widetilde{\cal B}^{\sigma}={\cal B}$,  the classical Bargmann transform.  In fact, 
let  $b$ and  $b^{\sigma}$ be functions  in one   complex variable and let 
$\widetilde{\cal B}: L^2({\bboard R}^r)\rightarrow {\cal F}({\bboard C}^r)$ and  
$\widetilde{\cal B}^{\sigma} : L^2({\bboard R}^r)\rightarrow {\cal F}({\bboard C}^r)$
 be the integral operators given by   : for $x\in S$, $\xi_z=zz^t \in \Gamma(\Xi)$,
$$(\widetilde{\cal B}f)(\xi_z)=\int_{S}b^{}(z,x)f(x)s(dx)$$
\noindent
and
$$(\widetilde{\cal B}^{\sigma}f)(\xi)=\int_{S}b^{\sigma}(z,x)f(x)s(dx).$$
\noindent
They  map  ${\cal C}^{\infty}(S)$ into ${\cal O}(\Gamma(\Xi))$, are  $O(r)$-equivariant, and map  ${\cal Y}_{k}({\bboard R}^{r})$ onto $\tilde{\cal O}_{k}({\bboard C}^r)$.
\vskip 4 pt
\noindent
The 'intertwining' relation for $\widetilde{\cal B}^{\sigma}$:

\vskip 6 pt

\centerline{$T(g')\widetilde{\cal B}^{\sigma}=\widetilde{\cal B}^{\sigma}R^{\sigma}(g)$}

\vskip 4 pt
\noindent
leads in particular to the  'intertwining' relations

\vskip 6 pt

\centerline{$dT(\widetilde E+\widetilde F)\circ\widetilde{\cal B}^{\sigma}=\widetilde{\cal B}^{\sigma}\circ dR^{\sigma}({\rm Ad}(g_0^{-1})(\widetilde E+\widetilde F))$,}

\vskip 4 pt
\noindent
and
\vskip 4 pt

\centerline{$dT(i(\widetilde E-\widetilde F))\circ\widetilde{\cal B}^{\sigma}=\widetilde{\cal B}^{\sigma}\circ dR^{\sigma}({\rm Ad}(g_0^{-1})(i(\widetilde E-\widetilde F))$.}

\vskip 4 pt
\noindent
But, since
\vskip 4 pt

\centerline{$\widetilde E={1\over 2}\pmatrix{0&0\cr
I_r&0}, \quad \widetilde F={1\over 2}\pmatrix{0&I_r\cr
0&0},$}

\vskip 4 pt
\noindent
then 
\vskip 4 pt

\centerline{$ \widetilde E+\widetilde F={1\over 2}\pmatrix{0&I_r\cr
I_r&0}$ and $i(\widetilde E-\widetilde F)={1\over 2}\pmatrix{0&-iI_r\cr
iI_r&0}$,}

\vskip 4 pt
\noindent
and
\vskip 4 pt

\centerline{${\rm Ad}(g_0^{-1})(\widetilde E+\widetilde F)={1\over 2}\pmatrix{0&-I_r\cr
-I_r&0}=-(\widetilde E+\widetilde F)$,}

\vskip 4 pt

\centerline{${\rm Ad}(g_0^{-1})(i(\widetilde E-\widetilde F))={1\over 2}\pmatrix{I_r&0\cr
0&-I_r}=:{1\over 2}\omega(I_r)$.}

\vskip 6 pt
\noindent
The 'intertwining' relations become
\vskip 4 pt

\centerline{$dT(\widetilde E)\circ\widetilde{\cal B}^{\sigma}+dT(\widetilde F)\circ\widetilde{\cal B}^{\sigma}=-\widetilde{\cal B}^{\sigma}\circ dR^{\sigma}(\widetilde E)-\widetilde{\cal B}^{\sigma}\circ dR^{\sigma}(\widetilde F)$,}

\vskip 4 pt
\noindent
and
\vskip 4 pt

\centerline{$idT(\widetilde E)\circ\widetilde{\cal B}^{\sigma}-idT(\widetilde F)\circ\widetilde{\cal B}^{\sigma}=\widetilde{\cal B}^{\sigma}\circ dR^{\sigma}({1\over 2}\omega(I_r)).$}

\vskip 6 pt
\noindent
It follows that 

\vskip 6 pt

\centerline{$2dT(\widetilde F)\circ\widetilde{\cal B}^{\sigma}=\widetilde{\cal B}^{\sigma} \circ\bigr({i\over 2} dR^{\sigma}(\omega(I_r))+dR^{\sigma}(-\widetilde E)+dR^{\sigma}(-\widetilde F)\bigl)$ \quad (*).}

\vskip 4 pt
\noindent

\vskip 8 pt
\noindent
Finally, using the integral form for the operator $\widetilde{\cal B}^{\sigma}$,  and the formulas

\vskip 4 pt

\centerline{$dR^{\sigma}(\omega(I_r)f(x)=dR^{\sigma}(\Bigr(\pmatrix{I_r&0\cr
0&-I_r}\Bigl)f(x)=-{r\over 2}f(x)+{\cal E}f(x)$ for $f\in L^2({\bboard R}^r),$}

\vskip 4 pt

\centerline{$dR^{\sigma}(\widetilde E)f(x)=-{i\over 4}\tau({\partial^2\over\partial x^2})f(x)  , \quad dR^{\sigma}(\widetilde F)f(x)=-{i\over 4}\tau(x^2)f(x)  $  for $f\in L^2({\bboard R}^r)$,}

\vskip 4 pt
\noindent
and
\vskip 4 pt

\centerline{$dT(\widetilde E)\phi(z)={i\over 4}\tau(z^2)\phi(z)$, \quad $dT(\widetilde F)\phi(z)={i\over 4}\tau({\partial^2\over\partial z^2})\phi(z)$  for $\phi\in {\cal F}({\bboard C}^r)$,}

\vskip 6 pt
\hfill
\eject
\noindent
one can deduce from (*)  that 
$$\eqalignno{&{i\over 2}\tau({\partial^2\over\partial z^2})\int_{{\bboard R}^r}b^{\sigma}(z,x)f(x)dx\cr
&=\int_{{\bboard R}^r}b^{\sigma}(z,x)\bigr(-{ir\over 4}+{i\over 2}\tau(x{\partial\over\partial x})+{i\over 4}\tau({\partial^2\over\partial x^2})+{i\over 4}\tau(x^2)\bigl)f(x)dx\cr
&=\int_{{\bboard R}^r}\bigr(-{ir\over 4}-{i\over 2}\tau(x{\partial\over\partial x})+{i\over 4}\tau({\partial^2\over\partial x^2})+{i\over 4}\tau(x^2)\bigl)b^{\sigma}(z,x)f(x)dx}$$

\vskip 4 pt

\vskip 4 pt
\noindent
which leads to the  following differential equation  for the   function $b^{\sigma}$:

\vskip 4 pt

\centerline{$-\tau({\partial^2\over\partial z^2})b^{\sigma}(z,x)=({r\over 2}+\tau(x{\partial\over\partial x})-{1\over 2}\tau({\partial^2\over\partial x^2})-{1\over 2}\tau(x^2)  \bigl)b^{\sigma}(z,x).$}

\vskip 4 pt
\noindent

\vskip 4 pt
\noindent
Observe that the solution of this equation is given by:
$$b^{\sigma}(z,x)=e^{-{1\over 2}(\tau(x^2)+\tau(z^2))+\sqrt 2\tau(zx)}.$$

\vskip 4 pt
\noindent
In fact, one has
\vskip 4 pt
$$\eqalign{{\partial^2\over\partial z_i^2}b(z,x)&={\partial\over\partial z_i}(-z_i+\sqrt 2 x_i)b^{\sigma}(z,x)\cr
&=(-1+(-z_i+\sqrt 2x_i)^2)b^{\sigma}(z,x)\cr
&=(-1+z_i^2+2x_i^2-2\sqrt 2x_iz_i)b^{\sigma}(z,x).}$$

\vskip 4 pt
\noindent
Then
\vskip 4 pt

\centerline{$-\tau({\partial^2\over\partial z^2})b^{\sigma}(z,x)=(r-\tau(z^2)-2\tau(x^2)+2\sqrt 2\tau(xz))b^{\sigma}(z,x).$}

\vskip 4 pt
\noindent
Similarily, one gets
\vskip 4 pt

\centerline{$-\tau({\partial^2\over\partial x^2})b^{\sigma}(z,x)=(r-\tau(x^2)-2\tau(z^2)+2\sqrt 2\tau(xz))b^{\sigma}(z,x).$} 

\vskip 6 pt
\noindent
On another part, 
\vskip 4 pt

\centerline{$x_i{\partial\over\partial x_i}b^{\sigma}(z,x)=-x_i^2+\sqrt 2x_iz_i$.}

\vskip 4 pt
\noindent
Then 
\vskip 4 pt

\centerline{$\tau(x{\partial\over\partial x})b^{\sigma}(z,x)=(-\tau(x^2)+\sqrt 2\tau(xz))b(z,x)$.}

\vskip 6 pt
\noindent
It follows that 
\vskip 4 pt
$$\eqalignno{&({r\over 2}+\tau(x{\partial\over\partial x})-{1\over 2}\tau({\partial^2\over\partial x^2})-{1\over 2}\tau(x^2)  \bigl)b^{\sigma}(z,x)\cr
&=({r\over 2}-\tau(x^2)+\sqrt 2\tau(zx)+{r\over 2}-{1\over 2}\tau(x^2)-\tau(z^2)+\sqrt 2\tau(zx)-{1\over 2}\tau(x^2))b^{\sigma}(z,x)\cr
&=(r-2\tau(x^2)-\tau(z^2)+2\sqrt 2\tau(zx))b^{\sigma}(z,x). \quad \qed}$$

\vskip 6 pt
\noindent

\bigskip
\centerline{\bf References}
\vskip 6pt
\noindent
\rm
\vskip 4 pt
\noindent
[A11] D. Achab (2011), \it Construction process for simple Lie algebras, \rm
J.  of Algebra, \bf 325,\rm 186-204.
\vskip 4 pt
\noindent
[AF12] D. Achab and J. Faraut (2012), \it Analysis of the Brylinski-Kostant model for minimal representations, Canad.J. of Math.,\bf 64, \rm721-754.
\vskip 4 pt
\noindent
[A12a]  D. Achab (2012), \it  Minimal representations of simple real Lie groups of non Hermitian type, \rm arxiv.
\vskip 4 pt
\noindent
[A12b]  D. Achab (2012), \it  Minimal representations of simple real Lie groups of  Hermitian type, \rm arxiv.
\vskip 4 pt
\noindent
[A16]  D. Achab  (2016), \it Analysis of the minimal representation of ${\rm SL}(n,{\bboard R})$, \rm arxiv.
\vskip 4 pt
\noindent
[AF16]  D. Achab and J. Faraut (2016), \it Analysis of the minimal representation of $O(m,n)$, \rm Preprint.
\vskip 4 pt
\noindent
[B97] R. Brylinski (1997), \it Quantization of the 
4-dimensional nilpotent orbit of ${\rm SL}(3,{\bboard R})$, \rm 
Canad.  Jou.l of Math.,\bf 49,\rm 916-943.
\vskip 4 pt
\noindent
[B98]  R. Brylinski (1998), \it Geometric quantization of real minimal nilpotent orbits, \rm Symp. Geom.,
Diff. Geom.  and Appl., \bf 9, \rm 5-58.
\vskip 4pt
\noindent
[BK94] R. Brylinski and B. Kostant (1994), \it 
Minimal representations, geometric quantization and unitarity, \rm 
Proc. Nat. Acad., \bf 91, \rm 6026-6029.
\vskip 4 pt
\noindent
[F15]  J. Faraut (2015), \it Analysis of the minimal representation of $O(n,n)$, \rm Preprint.
\vskip 4 pt
\noindent
[KV78] M. Kashiwara and M. Vergne (1978), \it  On the Segal-Shale-Weil Representations and Harmonic  Polynomials, \rm Invent. Math. \bf 44, \rm 1-47.
\vskip 4 pt
\noindent
[KO03a] T. Kobayashi and B. Orsted (2003), \it Analysis on the minimal representation of $O(p,q)$, I. Realization via conformal geometry, \rm Adv. Math., \bf 180, \rm  486-512.
\vskip 4 pt
\noindent
[KO03b] T. Kobayashi and B. Orsted (2003), \it Analysis on the minimal representation of $O(p,q)$, II. Branching laws, \rm Adv. Math.,\bf 180, \rm 513-550.
[KO03c] T. Kobayashi and B. Orsted (2003), \it Analysis on the minimal representation of $O(p,q)$, III.  Ultrahyperbolic equations on ${\bboard R}^{p-1,q-1}$, \rm Adv. Math., \bf 180, \rm 551-595.
\vskip 4 pt
\noindent
[KM11] T. Kobayashi and G. Mano (2011), \it The Schr\"odinger model for the minimal representation of the indefinite orthogonal group $O(p,q)$, \rm Memoir of the American Mathematical Society, Volume 213, Number 1000.
\vskip 4 pt
\noindent
[M] V.F. Molchanov (1970), \it Representations of pseudo-orthogonal groups associated with a cone, \rm Math. USSR Sbornik, \bf 10, \rm 333-347.

\end{document}

\end{document}